\documentclass{article}
\usepackage{amssymb,amsfonts,amsmath,amsthm,amscd}
\usepackage{makeidx}
\usepackage[english]{babel}
\usepackage{schemata}
\makeindex
\usepackage{tikz-cd}  
\usepackage{graphicx}
\usepackage{multicol}
\usepackage{tikz}

\usepackage{hyperref}
\input xy
\xyoption{all}
\usepackage[all]{xy}

\usepackage{fancyhdr}
\usepackage{makeidx}
\makeindex 
\usepackage{vmargin} 
\setmargins{3.3cm}       
{3cm}                        
{14.5cm}                      
{23cm}                    
{10pt}                           
{0.5cm}                           
{0pt}                             
{1cm}                           

\title{\textbf{About the space of continuous functions with open domain}}
\author{
	E. Ramírez\\
	\small  Mathematics School\\
	\small Industrial University of Santander\\
	\small Cl. 14 \# 6 - 07, Socorro, Santander, Colombia\\
	\small  e-mail: edwar5119@gmail.com\\
        \date{}
}
\theoremstyle{definition}
  \newtheorem{defn}{Definition}
  \newtheorem{thm}[defn]{Theorem}
  \newtheorem{obs}[defn]{Observation}
  
  \newtheorem{prop}[defn]{Proposition}
  \newtheorem{coro}[defn]{Corollary}
  \newtheorem{lema}[defn]{Lemma}

\begin{document}
\maketitle
\begin{abstract}
The classical compact-open topology is extended to the space of continuous functions with open domain using the Fell topology. Specifically, we will show that this extension makes the space completely metrizable. In particular, we will present a metric for the inverse semigroup of homeomorphisms of a locally compact, Hausdorff, and second-countable space.
\end{abstract}
	\textbf{2024 AMS 
subject to classification:} 
\section{Introduction} 
Functions with open domains appear in a variety of contexts in mathematics. Motivated by this phenomenon, Ahmed M. Abd-Allah and Ronald Brown set out to endow the space of continuous functions with open domains with an appropriate topology. In chronological order,
in \cite{AD}, they define the topology $\tau_{co}$, which 
turns the inverse semigroup of partial homeomorphisms between open subsets into a topological inverse semigroup, called $\Gamma(X)$ for $X$ locally compact Hausdorff space. More recently, in \cite{GX}, they added open sets after realizing that this topology is not $T_1$, showing that this addition makes $\Gamma(X)$ a Hausdorff topological inverse semigroup. Then, in \cite{Ston}, they show that if $X$ is a compact, Hausdorff, and totally disconnected space, then $(\Gamma(X),\tau_{hco})$ is completely metrizable.
In what follows, we will see how to define such a topology, we will see that the space of continuous functions with open domains called $C_{od}(X,Y)$ is a completely metrizable space when $X$ is  locally compact, Hausdorff and second countable space and $Y$ is a complete metric space. In particular, we will find a metric for $(\Gamma(X),\tau_{hco})$.

In \cite{AD}, the authors define the topology $\tau_{co}$ whose subbasic open sets are of the form
$$\langle K,V\rangle=\{f\in\Gamma(X):K\subseteq \mathrm{dom}(f)\text{ and }f(K)\subseteq V\}$$
where $K$ is compact and $V$ is open in $X$. 
\begin{obs}
The topology $\tau_{co}$ can be defined on the larger set of continuous functions whose domains are open subsets of $X$, to values in a topological space $Y$. This space was introduced in \cite{GX} and was denoted by $C_{od}(X,Y)$.
\end{obs}
Note that if $f \in \langle K,V \rangle$, any continuous function that coincides with $f$ in $K$ belongs to $\langle K, V \rangle$. This implies that in general, this topology is not $T_1$. Additionally, in \cite{AD}, the authors define the weak topology on $\Gamma(X)$, denoted by $\tau_{ico}$ and
given by the functions:
\begin{itemize}
    \item[] $\iota:\Gamma(X)\longrightarrow(\Gamma(X),\tau_{co})$ the inclusion map.
    \item[] $\eta:\Gamma(X)\longrightarrow(\Gamma(X),\tau_{co})$ the inversion map.
\end{itemize}
that is, to $\tau_{co}$ 
the subbasic open sets of the form $$\langle K,V\rangle^{-1}=\{f\in\Gamma(X):K\subseteq \mathrm{im}(f)\text{ and }f^{-1}(K)\subseteq V\}$$
are added. 

In \cite{AD}, the authors show that $(\Gamma(X), \tau_{ico})$ is a topological inverse semigroup. Recently, in \cite{GX}, the authors proved that this topology turns out to be the smallest one that makes $\Gamma(X)$ a topological inverse semigroup. Additionally, by utilizing the Fell topology [see Definition \ref{CLX}] on the hyperspace $CL(X)$ of closed sets including the empty set, the authors define the topology $\tau_{hco}$ by essentially adding to the weak topology of $\tau_{ico}$ the functions
\begin{itemize}
    \item[] $D:\Gamma(X)\longrightarrow(CL(X),\tau_{_{Fell}});\ f\longmapsto D(f)=X\setminus \mathrm{dom}(f)$
    \item[] $I:\Gamma(X)\longrightarrow(CL(X),\tau_{_{Fell}});\ f\longmapsto I(f)=X\setminus \mathrm{im}(f)$
\end{itemize}
This topology turns out to be the smallest Hausdorff topology that makes $\Gamma(X)$ a topological inverse semigroup \cite{GX}.
\section{$(\Gamma(X),\tau_{_{hco}})$ as metrizable space}
In \cite{Ston}, the authors show that $(\Gamma(X),\tau_{hco})$ is a metrizable space when $X$ is a compact Hausdorff space. In what follows we show than when $X$ is a locally compact and Hausdorff space, $(\Gamma(X),\tau_{hco})$ is a metrizable space. In the next section we will present an explicit metric for this space.
\begin{defn}\label{CLX}
$CL(X)$ is the family of closed subsets of $X$ including the empty set. This space has as subbasic open those of the form $(X\setminus K)^+=\{A\in CL(X):A\subseteq X\setminus K\}$ and $V^-=\{A\in CL(X):A\cap V\neq\emptyset\}$ where $K$ is compact and $V$ is non-empty open subset of $X$. This topology is called the Fell topology.

That is, to $(\Gamma(X),\tau_{ico})$ 
the subbasic open sets of the form $$D^{-1}((X\setminus K)^+),\ D^{-1}(V^-),\ I^{-1}((X\setminus K)^+)\text{ and }I^{-1}(V^-)$$
are added.
\end{defn} 
The following statement is straightforward.
\begin{prop}
    If $X$ is a Hausdorff locally compact and second countable space, then there exists a countable basis in $X$ of \textit{relatively compact} open sets, \textit{i.e.}, with compact closure.
\end{prop}
\begin{lema}\label{3333}
If $X$ is a Hausdorff locally compact and second countable space, then $(\Gamma(X),\tau_{co})$ is second countable space.
\end{lema}
\begin{proof}
Fix a relatively compact basis $\mathcal{B}$ closed under finite unions, and given $f \in \langle K,W \rangle$. For each $x\in K$, there exists $V_x\in \mathcal{B}$ such that $f(x)\in V_x\subseteq \overline{V_x}\subseteq W$. 
Due to the compactness of $f(K)$, there exist $x_1, \ldots, x_n \in K$ such that
$$f(K)\subseteq \bigcup \limits_{i=1}^nV_{x_i}=V_m\in\mathcal{B}$$
since $\mathcal{B}$ is closed under finite unions.
$$V_m\subseteq \overline{V_m}=\bigcup\limits_{i=1}^m \overline{V_{x_i}}\subseteq W
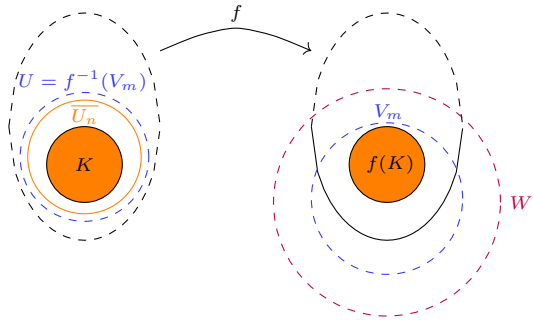
$$
\begin{center}
\begin{tikzpicture}
\draw[dashed] plot[domain=-3:-1] (\x,{1.5*sqrt(1-(\x+2)^2)}) plot[domain=-3:-1] (\x,{-1.5*sqrt(1-(\x+2)^2)}) plot[domain=1:3] (\x,{1.5*sqrt(1-(\x-2)^2)});
\draw plot[domain=1:3] (\x,{-1.5*sqrt(1-(\x-2)^2)});
\draw[dashed,blue!80] (-2,-0.4) circle (0.85)  (2,-0.95) circle (1) (1.7,0.2) node[anchor=west]{\scriptsize{$V_m$}}
(-3,0.6) node[anchor=west]{\scriptsize{$U=f^{-1}(V_m)$}};
\draw[fill=orange!100](-2,-0.5) circle (0.5) (2,-0.5) circle (0.5);
\draw[color=orange!100] (-2,-0.4) circle (0.75) (-2.3,0.16) node[anchor=west]{\scriptsize{$\overline{U_n}$}};
\draw[dashed, color=purple!100] (2,-1) circle (1.5) (3.5,-1) node[anchor=west]{\scriptsize{$W$}};
\draw (-0.2,1.5) node[anchor=west]{\scriptsize{$f$}} (-2.25,-0.5) node[anchor=west]{\scriptsize{$K$}} (1.6,-0.5) node[anchor=west]{\scriptsize{$f(K)$}};
\draw[->] plot[smooth] coordinates
{(-1,1)(-0.5,1.2)(0,1.3)(0.5,1.2)(1,1)};
\end{tikzpicture} 
\end{center}
Let $U=f^{-1}(V_m)$, $U$ be an open subset of $\mathrm{dom}(f)$, and therefore it is open in $X$. Since $K\subseteq U$, there exists $U_n\in \mathcal{B}$ 
such that $K\subseteq U_n\subseteq\overline{U_n}\subseteq U$. So $f$ belongs to $\langle \overline{U_n},V_m\rangle$.\\
If $g\in\langle \overline{U_n},V_m\rangle$,  
$K\subseteq \overline{U_n}\subseteq \mathrm{dom}(g)$ and $g(K)\subseteq g(\overline{U_n})\subseteq V_m\subseteq W$. So $\langle \overline{U_n}, V_m\rangle \subseteq \langle K,W\rangle$ and therefore 
$$\mathfrak{B}=\{\langle \overline{U},V\rangle:U,V\in \mathcal{B}\}$$
is a countable subbasis of $(\Gamma(X),\tau_{co})$.
\end{proof}
\begin{obs}\label{22}
In the previous proof only
the continuity  of $f$ was used, so $(C_{od}(X,Y),\tau_{co})$ is second countable space when $X$ and $Y$ are second countable spaces with $X$ locally compact Hausdorff space \cite[Proposition 2.1]{Mol}.
\end{obs} 
\begin{lema}\label{L2}
If $X$ is a locally compact Hausdorff space, then $(\Gamma(X),\tau_{hco})$ is regular.
\end{lema}
\begin{proof}
Fix a relatively compact basis $\mathcal{B}$ closed under finite unions.
Let $f\in \langle K,W\rangle$ with $K\neq\emptyset$ (since if $K=\emptyset$, then $\langle \emptyset,W\rangle=\Gamma(X)$).
If $x\in K$, there exists $V_x\in \mathcal{B}$ such that $f(x)\in V_x\subseteq\overline{V_x}\subseteq W$. Due to the compactness of $f(K)$, there exist $x_1,...,x_n\in K$ such that $$f(K)\subseteq \bigcup\limits_{i=1}^n V_{x_i}=V\subseteq\overline{V}=\bigcup\limits_{i=1}^n \overline{V_{x_i}}\subseteq W,\ V\in\mathcal{B}$$
Similarly, there exists $U \in \mathcal{B}$ such that $K\subseteq U\subseteq \overline{U}\subseteq f^{-1}(\overline{V})\subseteq \mathrm{dom}(f)$. This is $$f\in \langle K,V\rangle\cap D^{-1}((X\setminus\overline{U})^{+})=M$$
Let $g\in \overline{M}$ ($\overline{M}=cl_{\tau_{hco}}(M)$), there exists a red $(g_\lambda)$ in $M$ such that $\xymatrix{g_\lambda \ar@{->}[r]^{\tau_{hco}} & g}$. Let's see that $K\subseteq \mathrm{dom}(g)$. If we assume that $K\nsubseteq \mathrm{dom}(g)$, then $(X\setminus \mathrm{dom}(g))\cap K\neq \emptyset$, that is, $(X\setminus \mathrm{dom}(g))\cap U\neq \emptyset,\ (K\subseteq U)$ and therefore $D(g)=X\setminus \mathrm{dom}(g)\in U^-$. Since $\xymatrix{D(g_\lambda) \ar@{->}[r]^{\tau_{Fell}} & D(g)}$, there exists $\lambda_0$ such that $D(g_{_{\lambda_0}})=X\setminus dom(g_{_{\lambda_0}})\in U^-$, ie.
$$X\setminus dom(g_{_{\lambda_0}})\cap U\neq \emptyset $$
Which is a contradiction since $g_{_{\lambda_0}}\in D^{-1}((X\setminus \overline{U})^+)$, so $K\subseteq \mathrm{dom}(g)$. If we now assume that $g\notin \langle K,\overline{V}\rangle$, how $K\subseteq \mathrm{dom}(g)$, then $g(K)\nsubseteq \overline{V}$.
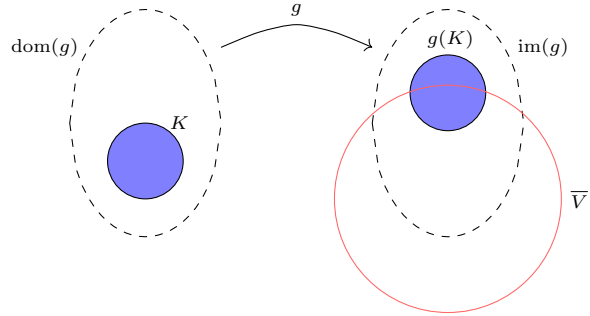
\begin{center}
\begin{tikzpicture}
\draw[dashed] plot[domain=-3:-1] (\x,{1.5*sqrt(1-(\x+2)^2)}) plot[domain=-3:-1] (\x,{-1.5*sqrt(1-(\x+2)^2)}) plot[domain=1:3] (\x,{1.5*sqrt(1-(\x-2)^2)}) plot[domain=1:3] (\x,{-1.5*sqrt(1-(\x-2)^2)});
\draw[fill=blue!50](-2,-0.5) circle (0.5) (2,0.4) circle (0.5);
\draw[color=red!60] (2,-1) circle (1.5);
\draw (-3.9,1) node[anchor=west]{\scriptsize{$\mathrm{dom}(g)$}} (2.8,1) node[anchor=west]{\scriptsize{$\mathrm{im}(g)$}}
(-1.8,0) node[anchor=west]{\scriptsize{$K$}} (1.6,1.1) node[anchor=west]{\scriptsize{$g(K)$}} (3.5,-1) node[anchor=west]{\scriptsize{$\overline{V}$}} (-0.2,1.5) node[anchor=west]{\scriptsize{$g$}};
\draw[->] plot[smooth] coordinates
{(-1,1)(-0.5,1.2)(0,1.3)(0.5,1.2)(1,1)};
\end{tikzpicture} 
\end{center}
There exists $a\in K$ such that $g(a)\in \mathrm{im}(g)\setminus \overline{V}$, that is, $g\in \langle \{a\},\mathrm{im}(g)\setminus\overline{V}\rangle$, then there exists $\lambda_0$ such that $g_{_{\lambda_0}}\in \langle \{a\},\mathrm{im}(g)\setminus\overline{V}\rangle$, thus $g_{_{\lambda_0}}(a)\notin\overline{V}$. But this is a contradiction because $g_{_{\lambda_0}}\in \langle K,V\rangle$. So $g \in \langle K,\overline{V}\rangle$ and it follows that $f \in M\subseteq\overline{M}\subseteq\langle K,\overline{V}\rangle\subseteq\langle K,W\rangle$.\\
If $f\in\langle K, W\rangle^{-1}$,  $f^{-1}\in\langle K, W\rangle$ And from the above, there exists $M$ $\tau_{hco}$-open in $\Gamma(X)$ such that $f^{-1}\in M\subseteq \overline{M}\subseteq\langle K,W\rangle$. Then $f=(f^{-1})^{-1}\in M^{-1}$ that is $\tau_{hco}$-open. Since inversion is $\tau_{hco}-\tau_{hco}$ homeomorphism, it follows that $$\left(\overline{M}\right)^{-1}=\overline{M^{-1}}$$
and therefore $f\in M^{-1}\subseteq\overline{M^{-1}}\subseteq\langle K, W\rangle^{-1}$.
Let $f\in D^{-1}(M)$, with $M\ \tau_{_{Fell}}-$open of $CL(X)$. As $X$ is locally compact Hausdorff space, by \cite[5.1.4 Corollary]{Beer} $(CL(X),\tau_{_{Fell}})$ is compact Hausdorff space, so it is regular, and therefore there exists an open subset $N$ of $CL(X)$ such that $D(f)\in N\subseteq \overline{N}^{Fell}\subseteq M$, then $$f\in D^{-1}(N)\subseteq \overline{D^{-1}(N)}^{hco}\subseteq D^{-1}\left(\overline{N}^{Fell}\right)\subseteq D^{-1}(M)$$ 
The case $f\in I^{-1}(M)$, with $M\ \tau_{_{Fell}}-$open of $CL(X)$ is analogous to the previous case.
\end{proof}
 \begin{obs}\label{22}
     Note that $\langle K,Y\rangle=D^{-1}((X\setminus K)^+)$ and $\langle K,X\rangle^{-1}=I^{-1}((X\setminus K)^+)$, for each $K$ compact subset of $X$. This means that the subsets added to $\tau_{ico}$ for generate the topology $\tau_{hco}$ are only subsets of the form $D^{-1}(V^-)$ and $I^{-1}(W^-)$ that is, neighborhoods of the function $\emptyset$. 
 \end{obs}
\begin{thm}
If $X$ is locally compact Hausdorff and second countable space, then $(\Gamma(X),\tau_{hco})$ is a separable metrizable space. 
\end{thm}
\begin{proof}
    From \cite[5.1.5 Theorem]{Beer} $(CL(X),\tau_{_{Fell}})$ is compact and metrizable, hence it is second countable. Let us prove that $(\Gamma(X),\tau_{hco})$ is second countable. Let $U$ be a $\tau_{hco}$-open set in $\Gamma(X)$ with $f\in U$, we have two cases:
    
    If $f=\emptyset$ then there exist $V_1,...,V_n,W_1,...W_r$ non-empty open subsets of $X$ such that (the basic neighborhoods of the empty function are of this form)
    $$\emptyset\in D^{-1}(V_1^-)\cap\cdots \cap D^{-1}(V_n^-)\cap I^{-1}(W_1^-)\cap\cdots \cap I^{-1}(W_r^-)\subseteq U$$
    Consider a countable basis $\mathcal{B}$ relatively compact and closed under finite unions, and we can assume that such $V_i$'s and $W_j$'s are elements of $\mathcal{B}$.
    
    If $f\neq\emptyset$, using $\mathcal{B}$ and \cite[5.1.5 Theorem]{Beer} taking the countable subbasis of $(CL(X),\tau_{_{Fell}})$ 
$$\{(X\setminus \overline{V})^+,\ V^-:V\in\mathcal{B}\}$$ and Lemma \ref{3333}, there exists $R_i,S_i,T_j,U_j,V_k,W_l\in\mathcal{B}$ such that 
    $$\begin{array}{lcllll}
       f & \in & \langle\overline{R_1},S_1\rangle&\cap\cdots\cap&\langle\overline{R_m},S_m\rangle&\cap\\
       &&\langle\overline{T_1},U_1\rangle^{-1}&\cap\cdots\cap&\langle\overline{T_n},U_n\rangle^{-1}&\cap\\
       &&D^{-1}(V_1^-)&\cap\cdots\cap&D^{-1}(V_r^-)&\cap\\
   &&I^{-1}(W_1^-)&\cap\cdots\cap&I^{-1}(W_s^-)&\subseteq U.
    \end{array}
    $$
   So for the Observation  \ref{22} we can consider the \textit{\textbf{countable subbasis}} $$\mathfrak{B}=\{\langle \overline{U},V\rangle,\ \langle \overline{U},V\rangle^{-1},\ D^{-1}(V^-),\ I^{-1}(W^-):U,V,W\in \mathcal{B}\}$$
    By \cite[Proposition 3.10]{GX} $(\Gamma(X),\tau_{hco})$ is Hausdorff space (and therefore $T_1$), by Lemma \ref{L2} and using Urysohn's metrization theorem, we conclude that $(\Gamma(X),\tau_{hco})$ is a metrizable space.
\end{proof}
\section{About the space $(C_{od}(X,Y),\tau_{_{\iota,D}})$}
We equip the space $C_{od}(X,Y)$ with the weak topology given by the functions $\iota$ and $D$ where
\begin{itemize}
    \item[] $\iota:C_{od}(X,Y)\longrightarrow (C_{od}(X,Y),\tau_{co})$ the inclusion;
    \item[] $D:C_{od}(X,Y)\longrightarrow (CL(X),\tau_{_{Fell}}),\ D(f)=X\setminus \mathrm{dom}(f).$
\end{itemize}
called $\tau_{_{\iota,D}}$. In what follows, $\mathbb{K}(X)$ will denote the hyperspace of non-empty compact subsets of a topological space $X$. 
\begin{obs}\label{2n}
Note that if $\mathcal{B}_1$ and $\mathcal{B}_2$ are countable basis and closed under finite unions of 
$X$ and $Y$ respectively, and $\mathcal{B}_1$ is relatively compact, then the set
$$\mathfrak{C}=\{\langle \overline{U},V\rangle,\ D^{-1}(U^-):U\in \mathcal{B}_1 \ \text{and}\ V\in \mathcal{B}_2\}$$ is a countable subbasis for $(C_{od}(X,Y),\tau_{_{\iota,D}})$. 
\end{obs}
\subsection{Convergence in $C_{od}(X,Y)$}
Let $X$ be a topological space and $(Y,d)$ be a metric space, $f,g\in C_{od}(X,Y)$ with $\mathrm{dom}(f)\cap \mathrm{dom}(g)\neq\emptyset$. Define
$$d_K(f,g):=\sup\limits_{x\in K}d(f(x),g(x)),\ \  K\in\mathbb{K}(\mathrm{dom}(f)\cap \mathrm{dom}(g))$$
\begin{defn}
    Let $f\in C_{od}(X,Y)\setminus\{\emptyset\}=C_{od}^\star(X,Y)$, $K$ non-empty compact subset of $\mathrm{dom}(f)$ and $\epsilon > 0$, we define 
        $$B_K(f,\epsilon):=\{g\in C_{od}^\star(X,Y):K\in \mathbb{K}(\mathrm{dom}(g))\text{ and }d_K(f,g)< \epsilon\}
        $$
\end{defn}
\begin{lema}
    $\mathcal{B}=\{B_K(f,\epsilon):\, f\in C_{od}^\star(X,Y),\, K\in\mathbb{K}(\mathrm{dom}(f))\text{ and }\epsilon>0\}\cup\{C_{od}(X,Y)\}$ is a basis for a topology.
\end{lema}
\begin{proof}
    (i)$\bigcup\limits_{B\in \mathcal{B}}B=C_{od}(X,Y)$. \\
    (ii) Let $f\in B_K(g,\epsilon_1)\cap B_L(h,\epsilon_2)$ ($g,h\neq \emptyset$, $K\subseteq \mathrm{dom}(f)$ and $L\subseteq \mathrm{dom}(g)$). Consider $\delta_1=\epsilon_1-d_K(f,g)$ and $\delta_2=\epsilon_2-d_L(f,h)$\\
    \textbf{Statement 1.} $B_K(f,\delta_1)\subseteq B_K(g,\epsilon_1)$. Indeed, if $p\in B_K(f,\delta_1)$, $K\subseteq dom(p)$ and $d_K(f,p)<\delta_1=\epsilon_1-d_K(f,g)$ ie. $d_K(p,g)\leq d_K(f,p)+d_K(f,g)<\epsilon_1$, then $p\in B_K(g,\epsilon_1)$.\\
    \textbf{Statement 2.} $B_L(f,\delta_2)\subseteq B_L(h,\epsilon_2)$ The proof is analogous to the previous statement.\vspace{2mm}\\
    Take $\delta=\min\{\delta_1 ,\delta_2\}$. $f\in B_{K\cup L}(f,\delta)$ and if $p\in B_{K\cup L}(f,\delta)$, then $d_K(f,p)<\delta_1$ and $d_L(f,p)<\delta_2$ ie. $p\in B_{K}(f,\delta_1)\cap B_{L}(f,\delta_2)\subseteq B_K(g,\epsilon_1)\cap B_L(h,\epsilon_2)$.
\end{proof}
Let $\tau_{cc}$ be the topology generated by $\mathcal{B}$, and we'll call it \textbf{\textit{the topology of compact convergence.}}
\begin{prop}
   If $X$ is locally compact Hausdorff and second countable space, and $(Y,d)$ be a metric space. Then $\tau_{co}$ and $\tau_{cc}$ matches on $C_{od}(X,Y)$.
\end{prop}
\begin{proof}
    $\tau_{co}\subseteq \tau_{cc}$. Indeed, let $f\in \langle K,V \rangle$, ie. $K\subseteq \mathrm{dom}(f)$ and $f(K)\subseteq V$. Then there exists $\epsilon >0$ such that the cloud $N(\epsilon, f(K))\subseteq V$. If $g\in B_K(f,\epsilon)$, then $K\subseteq \mathrm{dom}(g)$ and $d_K(f,g)<\epsilon$, that is, $d(f(x),g(x))<\epsilon,\ \forall x\in K$. If $x\in K$, $d(g(x),f(K))\leq d(f(x),g(x))<\epsilon$. So $g(x)\in N(\epsilon,f(K))$ ie. $g(K)\subseteq N(\epsilon,f(K))\subseteq V$ and therefore $$B_K(f,\epsilon)\subseteq \langle K,V \rangle$$
    $\tau_{cc}\subseteq \tau_{co}$. Indeed, let $f\in C_{od}(X,Y)$, $K$ is
    non-empty compact subset of $\mathrm{dom}(f)$, and $\epsilon>0$. If $x\in K$, there exists $V_x$ open neighborhood of $f(x)$ with diameter $<\epsilon$. By continuity of $f$ and local compactness of $X$, there exists an open neighborhood of $x$ and relatively compact $U_x$ with $\overline{U_x}\subseteq \mathrm{dom}(f)$ such that $f(\overline{U_x})\subseteq V_x$, ie. $f\in \langle \overline{U_x},V_x \rangle,\ \forall x\in K.$\\
    By compactness of $K$, there exists $x_1,...,x_n \in K$ such that $K\subseteq  \bigcup\limits_{i=1}^nU_{x_i}$,  with $K\cap U_{x_i}\neq\emptyset$ for each $i\in\{1,...,n\}$.
    
    Let $K_i=\overline{U_{x_i}}\cap K$, since $X$ is compactly generated, $K_i$ is closed in $K$ and therefore compact. Finally, let's see that
    $$W=\bigcap\limits_{i=1}^n\langle K_i,V_{x_i}\rangle\subseteq B_K(f,\epsilon)$$
    Let $g\in W$ and $x\in K$. There exists $i\in\{1,...,n\}$ such that $x\in U_{x_i}$, then $x\in K_i$ and as $g\in \langle K_i,V_{x_i}\rangle$, $g(x)\in V_{x_i}$. Since $f\in\langle K_i,V_{x_i}\rangle$ and $x\in K_i$, $f(x)\in V_{x_i}$ and as $diam(V_{x_i})<\epsilon$, $d(f(x),g(x))<\epsilon$. Therefore $g\in B_K(f,\epsilon)$.
\end{proof}
This is \textit{analogous} to what happens when we define the compact-open topology on $Y^X$ and restrict to $C(X,Y)$ where $Y$ is a metric space \cite[Theorem 46.8]{Mun}.
\begin{lema}\label{coro}
    Let $X$ be a locally compact Hausdorff space and $Y$ be a metric space. If $(f_\lambda)$ is a net in $C_{od}(X,Y)$ $\tau_{cc}$-converging to a function $f$, then $f$ is continuous.
\end{lema}
\begin{proof}
    If $f\neq\emptyset$ , let $K\in\mathbb{K}(\mathrm{dom}(f))$. Since 
    $\xymatrix{f_\lambda \ar[r]^{\tau_{cc}} & f}$, we have that $K\subseteq dom(f_\lambda)$ eventually, and
    $\xymatrix{f_n|_{_{K}} \ar[r] & f|_{_{K}}}$ uniformly. So  $f|_K$ is continuous, since $X$ is a $k-$space the result follows from \cite[43.10 Lemma]{Will}.
\end{proof}
In particular, if $(f_n)$ is a sequence of \textit{holomorphic functions} defined on open domains $\tau_{cc}$-converging to $f\in C_{od}(\mathbb{C})$ then $f$ is holomorphic. This implies that the set of holomorphic functions defined on open subsets of $\mathbb{C}$ is $\tau_{cc}$-closed.
\subsection{$(C_{od}(X,Y), \tau_{_{\iota,D}})$ as metric space}
In what follows, $X$ will be a  locally compact Hausdorff and second countable space, and $(Y,d)$ be a metric space with $d \leq 1$. Also, $\mathbb{N}$ will be the set of positive integers.
Note that $$\bigcap\limits_{U\in \tau\setminus\{\emptyset\}}D^{-1}(U^-)=\{\emptyset\}$$
Intuitively, this indicates that if $U$ is a non-zero open subset of $X$, two functions are close to $\emptyset$ (for $U$) if $f,g\in D^{-1}(U^-)$. Let's see that it is possible to define an explicit metric on $(C_{od}(X,Y),\tau_{_{\iota,D}})$.\newline 

Let $\mathcal{B}=\{U_n\}_{n\in \mathbb{N}}$ a basis of $X\ (U_0=\emptyset)$, and for each $m,n\in \mathbb{N}$, lets $K_{mn}$ non-empty compact subsets of $X$ such that 
\begin{equation*}
    K_{mn}\subseteq int(K_{(m+1)n}) \text{ and }\bigcup\limits_{m\in \mathbb{N}}K_{mn}=U_n,\ \forall m,n\in \mathbb{N}
\end{equation*}
We define
$$\beta_{mn}(f,g)=\left\{\begin{array}{ll}
    0,  & f,g\in D^{-1}(int(K_{(m+1)n})^-);\\
    d_{K_{mn}}(f,g), & f,g\notin D^{-1}(int(K_{(m+1)n})^-);\\
    1, & \text{other cases.}
\end{array}\right.$$
$\beta_{mn}$ 
is a pseudometric on $C_{od}(X,Y)$. Indeed, given $f,g,h\in C_{od}(X,Y)$, Let's call $$L=D^{-1}(int(K_{(m+1)n})^-)$$ 
we have that

$$\beta_{mn}(f,g)\left\{\begin{array}{l}
    \underset{f,g\in L}{=}  0\leq \beta_{mn}(f,h)+\beta_{mn}(h,g)\\\\
    \underset{f,g\notin L}{=} d_{K_{mn}}(f,g)\left\{\begin{array}{l}
        \underset{h\in L}{\leq} 1+1=\beta_{mn}(f,h)+\beta_{mn}(h,g)\\
        \underset{h\notin L}{\leq}  d_{K_{mn}}(f,h)+d_{K_{mn}}(h,g)=\beta_{mn}(f,h)+\beta_{mn}(h,g)
    \end{array}\right.\\\\
    \underset{f\in L\text{ and } g\notin L}{=} 1 \left\{ \begin{array}{l}
         \underset{h\in L}{=} 0+1=\beta_{mn}(f,h)+\beta_{mn}(h,g)\\
         \underset{h\notin L}{\leq} 1+\beta_{mn}(h,g)=\beta_{mn}(f,h)+\beta_{mn}(h,g)
    \end{array}\right.\\\\
    \underset{f\notin L\text{ and } g\in L}{=} 1\left\{ \begin{array}{l}
         \underset{h\in L}{=} 1+0=\beta_{mn}(f,h)+\beta_{mn}(h,g)\\
         \underset{h\notin L}{\leq} \beta_{mn}(f,h)+1=\beta_{mn}(f,h)+\beta_{mn}(h,g)
    \end{array}\right.
\end{array}\right.$$

Let $\beta_n(f,g)=\sum\limits_{m=1}^\infty 2^{-m}\beta_{mn}(f,g)$, and finally
$$\beta(f,g)=\sum\limits_{n=1}^\infty 2^{-n}\beta_n(f,g)=\sum\limits_{n=1}^\infty\sum\limits_{m=1}^\infty 2^{-(m+n)}\beta_{mn}(f,g)$$
$\beta$ is a metric on $C_{od}(X,Y)$. Indeed, if $f,g\in C_{od}(X,Y)$ are such that $\beta(f,g)=0$ or equivalently $\beta_{mn}(f,g)=0,\text{ for each }m,n\in \mathbb{N}$. 

If $f=\emptyset$, since $\emptyset\in D^{-1}(int(K_{mn})^-),\text{ for each } m,n\in\mathbb{N}$ and $\beta_{mn}(\emptyset,g)=0$ for each $m,n\in \mathbb{N}$, by definition of $\beta_{mn}$ necessarily 
$$g\in\bigcap\limits_{n\in\mathbb{N}}\bigcap\limits_{m \in\mathbb{N}}D^{-1}(int(K_{mn})^-)\subseteq\bigcap\limits_{n\in\mathbb{N}}D^{-1}(U_n^-)=\bigcap\limits_{U\in \tau\setminus\{\emptyset\}}D^{-1}(U^-)=\{\emptyset\}$$
therefore, $g=\emptyset$. If $g=\emptyset$, analogously $f=\emptyset$.

If $f,g\neq \emptyset$, let $x\in \mathrm{dom}(f)$. There exists $m,n\in \mathbb{N}$ such that $x\in K_{mn}\subseteq U_n\subseteq \mathrm{dom}(f)$, since $f\notin D^{-1}(int(K_{(m+1)n})^-)$ and $\beta_{mn}(f,g)=0$, by definition of the pseudometric $\beta_{mn}$, necessarily $g\notin D^{-1}(int(K_{(m+1)n})^-)$ ($\textit{i.e.}\ x\in \mathrm{dom}(g)$) and in particular,
$$d(f(x),g(x))\leq d_{K_{mn}}(f,g)=\beta_{mn}(f,g)=0,$$
then $f(x)=g(x)$. Análogously, if $x\in \mathrm{dom}(g)$, then $x\in \mathrm{dom}(f)$ and $f(x)=g(x)$. Therefore $f=g$.
\begin{thm}
    If $\mathcal{B}$ is closed under finite unions, then $\beta$ is a metric for $(C_{od}(X,Y),\tau_{_{\iota,D}})$.
\end{thm}
\begin{proof}
    If $\xymatrix{f_\lambda \ar[r]^{\tau_{\iota,D}} & f}$, given $\epsilon>0$ and let $N\in\mathbb{N}$ such that $2^{-N}<\frac{\epsilon}{2}$. If $n\leq N$, We have two cases:\\
    \textbf{Case 1.} $f\in D^{-1}(U_n^-).$  When $f\in  D^{-1}(int(K_{1n})^-)$ there exist $\sigma$ such that $f_\lambda\in D^{-1}(int(K_{1n})^-)$, for each $\lambda>\sigma.$ Since $D^{-1}(int(K_{1n})^-)\subseteq D^{-1}(int(K_{mn})^-),\ \forall m\geq 1$ 
    We have that $\beta_{mn}(f,f_\lambda)=0<\frac{\epsilon\,2^{n-1}}{N}\cdot 2^{-m}$, for each $m\geq 1$ and each $\lambda> \sigma$. So $$\beta_n(f,f_\lambda)=0<\frac{\epsilon\,2^{n-1}}{N},\text{ for each }\lambda>\sigma$$
    We can assume that $f\notin D^{-1}(int(K_{1n}))$ we can two news subcases.\vspace{2mm}\\
    \textbf{Case 1.1.} When $f\in D^{-1}(int(K_{2n}))$. There exist $\sigma $ such that $f_\lambda\in D^{-1}(int(K_{2n})^-)$, for each $\lambda>\sigma$, and we have that $\beta_{1n}(f,f_\lambda)=0<\frac{\epsilon\,2^{n-1}}{N}\cdot 2^{-1}$ for each $\lambda> \sigma$. Given that $D^{-1}(int(K_{2n})^-)\subseteq D^{-1}(int(K_{mn})^-),\ \forall m\geq 2$ we have that $\beta_{mn}(f,f_\lambda)=0<\frac{\epsilon\,2^{n-1}}{N}\cdot 2^{-m}$. So
    $$\beta_n(f,f_\lambda)=0<\frac{\epsilon\,2^{n-1}}{N},\text{ for each }\lambda>\sigma$$   
    \textbf{Case 1.2} When $f\notin D^{-1}(int(K_{2n}))$.

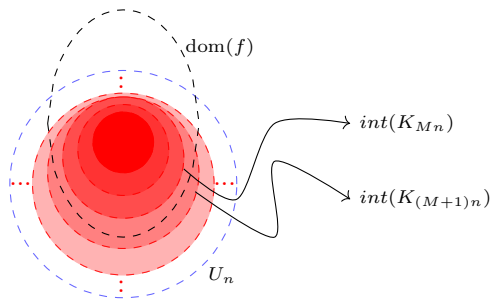
\begin{center}
\begin{tikzpicture}
\draw[dashed, color=blue!60] (2,-0.8) circle (1.5);
\filldraw[color=red!30] (2,-0.8) circle (1.2cm);
\filldraw[color=red!50] (2,-0.65) circle (1cm);
\filldraw[color=red!70] (2,-0.45) circle (0.8cm);
\filldraw[color=red!80] (2,-0.35)  circle (0.6cm) ;
\filldraw[color=red!100] (2,-0.25) circle (0.4cm);
\draw[dashed, color=red!90] (2,-0.8) circle (1.2cm);
\draw[dashed, color=red!95] (2,-0.65) circle (1cm);
\draw[dashed, color=red!96] (2,-0.45) circle (0.8cm);
\draw[dashed, color=red!98] (2,-0.35)  circle (0.6cm) ;
\draw[dashed, color=red!100] (2,-0.25) circle (0.4cm);
\draw[dashed] plot[domain=1:3] (\x,{1.5*sqrt(1-(\x-2)^2)}) plot[domain=1:3] (\x,{-1.5*sqrt(1-(\x-2)^2)});
\draw[->] plot[smooth] coordinates {(2.95,-0.9)(4,-1.5)(4,-0.5)(5,-1)};
\draw[->] plot[smooth] coordinates {(2.8,-0.6)(3.5,-1)(4,0)(5,0)};
\draw (5,0) node[anchor=west]{\scriptsize{$int(K_{Mn})$}} (5,-1) node[anchor=west]{\scriptsize{$int(K_{(M+1)n})$}} (2.7,1) node[anchor=west]{\scriptsize{$\mathrm{dom}(f)$}} (3,-2) 
node[anchor=west]{\scriptsize{$U_n$}} (8,1) node[anchor=west]{\scriptsize{There exists $M\geq 2$ such that}} (8,0.6) node[anchor=west]{\scriptsize{$f\in D^{-1}(int(K_{(M+1)n})^-)$ and }} (8,0.2) node[anchor=west]{\scriptsize{$f\notin D^{-1}(int(K_{Mn})^-)$.}};
\filldraw[color=red] (1.8,0.55) node[anchor=west]{:} (1.8,-2.15) node[anchor=west]{:} (0.37,-0.8) node[anchor=west]{...} (3.07,-0.8) node[anchor=west]{...};
\end{tikzpicture} 
\end{center}
On one hand, there exists $\sigma$ such that $f_\lambda\in D^{-1}(int(K_{(M+1)n})^-)$, $\forall \lambda >\sigma$. Given that $D^{-1}(int(K_{(M+1)n})^-)\subseteq D^{-1}(int(K_{mn})^-),\ \forall m\geq M+1$, it follows that $\beta_{mn}(f_\lambda,f
)=0,\ \forall \lambda >\sigma,\ \forall m\geq M$. We have that
$$\beta_n(f_\lambda,f)=\sum\limits_{m=1}^{M-1}2^{-m}\beta_{mn}(f_\lambda,f),\ \forall \lambda>\sigma.$$
On the other hand, $K_{1n}\subseteq K_{2n}\subseteq\cdots\subseteq K_{(M-1)n}\subseteq int(K_{Mn})\subseteq \mathrm{dom}(f)$, and since $\xymatrix{f_\lambda \ar[r]^{\tau_{cc}} & f}$, for each $m\in\{1,...,M-1\}$, there exists $\lambda_m$ such that
$$f_\lambda\in B_{K_{mn}}\left(f,\frac{\epsilon\,2^{m+n-1}}{MN}\right),\ \forall \lambda >\lambda_m,\,\forall m\in\{1,...,M-1\}$$
Let $\gamma_n>\lambda_1,...,\lambda_{M-1},\sigma$. Thus, we have that for each $\gamma >\gamma_n$
\begin{equation*}
    \begin{split}
    \beta_n(f_\lambda,f) & = \sum\limits_{m=1}^{M-1}2^{-m}\beta_{mn}(f_\lambda,f)\\
    & =\sum\limits_{m=1}^{M-1}2^{-m}d_{K_{mn}}(f_\lambda,f)\\
    & \leq\sum\limits_{m=1}^{M-1}2^{-m}\frac{\epsilon\,2^{m+n-1}}{MN}\\
    & =\epsilon\frac{M-1}{MN}2^{n-1}\  <\ \frac{\epsilon\,2^{n-1}}{N}
    \end{split}
\end{equation*}
\textbf{Case 2.} $f\notin D^{-1}(U_n^-)$, \textit{i.e.} $U_n\subseteq \mathrm{dom}(f)$.\\
Fix $M\in \mathbb{N}$ such that $2^{-M}< \frac{\epsilon\, 2^{n-2}}{N}$.
For each
$m\in\{1,...,M\}$, since $\xymatrix{f_\lambda \ar[r]^{\tau_{cc}} & f,}$there exists $\lambda_m$ such that $f_\lambda\in B_{K_{mn}}\left(f,\frac{\epsilon\, 2^{m+n-2}}{MN}\right),\ \forall\lambda>\lambda_m.$ Let $\gamma_n>\lambda_1,...,\lambda_m$. Thus, we have that for each $\lambda>\gamma_n$
\begin{equation*}
\begin{split}    
\beta_n(f_\lambda,f) & = \sum\limits_{m=1}^M 2^{-m}\beta_{mn}(f_\lambda,f)+\sum\limits_{m=M+1}^\infty 2^{-m}\beta_{mn}(f_\lambda,f) \\
    & \leq \sum\limits_{m=1}^M 2^{-m}\beta_{mn}(f_\lambda,f) +2^{-M}\\
    & < \sum\limits_{m=1}^M 2^{-m}\frac{\epsilon\, 2^{m+n-2}}{MN}+\frac{\epsilon\, 2^{n-2}}{N}\\
    &=\frac{\epsilon\ 2^{n-2}}{N}+\frac{\epsilon\ 2^{n-2}}{N}=\frac{\epsilon\ 2^{n-1}}{N}
\end{split}
\end{equation*}
Let $\gamma>\gamma_1,...,\gamma_n$. Thus, we have that for each $\lambda>\gamma$
\begin{equation*}
    \begin{split}
        \beta(f_\lambda,f) & = \sum\limits_{n=1}^\infty 2^{-n}\beta_{n}(f_\lambda,f)\\
        &= \sum\limits_{n=1}^N 2^{-n}\beta_{n}(f_\lambda,f)+\sum\limits_{n=N+1}^\infty 2^{-n}\beta_{n}(f_\lambda,f)\\
        &\leq \sum\limits_{n=1}^N 2^{-n}\frac{\epsilon\ 2^{n-1}}{N}+2^{-N}\\
        &<\frac{\epsilon}{2}+\frac{\epsilon}{2}=\epsilon
    \end{split}
\end{equation*}
That is, $\xymatrix{f_\lambda \ar[r]^{\beta} & f}$.\\ Conversely, if $\xymatrix{f_\lambda \ar[r]^{\beta} & f}$, then $\xymatrix{f_\lambda \ar[r]^{\tau_{\iota,D}} & f.}$  Indeed, lets $K\in\mathbb{K}(\mathrm{dom}(f))$ and $\epsilon>0$. For each $x\in K$, there exist $U_x\in \mathcal{B}$ such that $x\in U_x\subseteq \mathrm{dom}(f)$. By compactness of $K$, there exists $x_1,...,x_l\in K$ such that $$K\subseteq\bigcup\limits_{i=1}^l U_{x_i}=U_N\in\mathcal{B}\    (\mathcal{B}\ \text{is closed under finite unions).}$$
Since $K_{mN}\subseteq int(K_{(m+1)N}),\forall m\in \mathbb{N},$ there exists $M\in \mathbb{N}$ such that $K\subseteq K_{MN}\subseteq U_N\subseteq \mathrm{dom}(f)$. Due to the convergence in the hypothesis, there exists $\lambda_0$ such that $$\beta(f_\lambda,f)<\epsilon\cdot2^{-(M+N)},\text{ for each }\lambda>\lambda_0$$ 
So
\begin{equation*}
    \begin{split}
        d_K(f_\lambda,f) & \leq d_{K_{MN}}(f_\lambda,f) \\
        & =2^{M+N}\cdot 2^{-(M+N)}\beta_{MN}(f_\lambda,f)\\
        & \leq 2^{M+N}\beta(f_\lambda,f)\\
        & < 2^{M+N}\frac{\epsilon}{2^{M+N}}=\epsilon
    \end{split}
\end{equation*}
Therefore $\xymatrix{f_\lambda \ar[r]^{\tau_{cc}} & f}.$\\
Let $U$ nonempty open set in $X$, $f\in D^{-1}(U^-)$ and fix $0<\epsilon<1$.
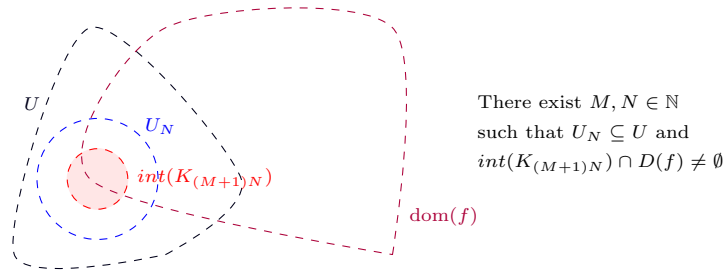
\begin{center}
    \begin{tikzpicture}
        \filldraw[color=red!10] (-0.9,0) circle (0.4);
        \draw[dashed, color=red] (-0.9,0) circle (0.4);
        \draw[dashed, color=blue!10!black] plot[smooth] coordinates
        {(0,-1)(1,0)(-1,2)(-2,-1)(0,-1)};
        \draw[color=blue!10!black] (-2,1) node[anchor=west]{\scriptsize{$U$}};
        \draw[color=red] (-0.5,0) node[anchor=west]{\scriptsize{$int(K_{(M+1)N})$}};
        \draw[dashed, color=purple!90!black] plot[smooth] coordinates 
        {(3,-1)(3,2)(0,2)(-1,0)(3,-1)};
        \draw[color=purple!90!black] (3.1,-0.5) node[anchor=west]{\scriptsize{$\mathrm{dom}(f)$}};
        \draw[dashed, color=blue] (-0.9,0) circle (0.8);
        \draw[color=blue] (-0.4,0.7) node[anchor=west]{\scriptsize{$U_N$}};
        \draw (4,1) node[anchor=west]{\scriptsize{There exist $M,N\in\mathbb{N}$}} (4,0.6) node[anchor=west]{\scriptsize{such that $U_N\subseteq U$ and}} (4,0.2) node[anchor=west]{\scriptsize{$int(K_{(M+1)N})\cap D(f)\neq \emptyset$}};
    \end{tikzpicture}
\end{center}
Due to $\beta$-convergence, there exists $\lambda_0$ such that $$\beta(f_\lambda,f)<\frac{\epsilon}{2^{M+N}},\ \forall \lambda>\lambda_0.$$
Thus, for each $\lambda>\lambda_0$ we have that
\begin{equation*}
    \begin{split}
    \beta_{MN}(f_\lambda,f) & =2^{M+N}\cdot 2^{-(M+N)}\beta_{MN}(f_\lambda,f)\\
    &\leq2^{M+N}\beta(f_\lambda,f)\\
    &<2^{M+N}\frac{\epsilon}{2^{M+N}}=\epsilon<1
    \end{split}
\end{equation*}
since $f\in D^{-1}(int(K_{(M+1)N})^-)$, by definition of $\beta_{MN}$, necessarily $f_\lambda\in D^{-1}(int(K_{(M+1)N})^-),$ for each $\lambda>\lambda_0$. Thus
$$\emptyset\neq int(K_{(M+1)N})\cap D(f_\lambda)\subseteq U_N\cap D(f_\lambda)\subseteq U\cap D(f_\lambda),\ \forall \lambda>\lambda_0.$$
Therefore $\xymatrix{f_\lambda \ar[r]^{\tau_{\iota,D}} & f}.$
\end{proof}
In particular, if $(X,d)$ is a locally compact and second countable metric space ($d\leq 1$), for each $f,g\in\Gamma(X)$ $$d_\gamma(f,g)=\beta(f,g)+\beta(f^{-1},g^{-1})$$
define a \textit{metric} for $(\Gamma(X),\tau_{hco})$.\vspace{2mm}\\
Also, using \cite[Proposition 3.14]{GX} we can define a metric for the Fell topology on $CL(X)$ as follows, for each $A,B\in CL(X)$
$$d_{_{Fell}}(A,B):=\frac{1}{2}d_\gamma(\textbf{Id}_{X\setminus A},\textbf{Id}_{X\setminus B})=\sum\limits_{n=1}^\infty\sum\limits_{m=1}^\infty 2^{-(m+n)}\left\{
\begin{array}{ll}
    0, & A,B\in int(K_{(m+1)n})^-\ \text{ or }\\
    & A,B\notin int(K_{(m+1)n})^-;\\
    &\\
    1, & \text{other cases.}
\end{array}\right.$$

\begin{obs}\label{aaa}
Note that $d_{_{Fell}}(D(f),D(g))\leq \beta(f,g)$, for each $f,g\in C_{od}(X,Y)$.
Thus, every $\beta$-Cauchy net is such that the net of the complements of their domains are $d_{_{Fell}}$-Cauchy. 
\end{obs}
\subsection{$(C_{od}(X,Y),\beta)$ as Polish space}
In this section, we will see that if $X$ is a locally compact, Hausdorff and second countable space, and $(Y,d)$ is a complete metric space ($d\leq 1$), then $(C_{od}(X,Y),\beta)$ is a complete metric space. In particular, \textit{the set of holomorphic functions from open domain in $\mathbb{C}$ will be a closed subspace and thus Polish with the metric $\beta$.}
\begin{defn}
    Let $X$ be a locally compact, Hausdorff and second countable space, and $Y$ be a complete metric space.
    A net $(f_\lambda)\subseteq C_{od}(X,Y)$ is called $\gamma-Cauchy$ if $(D(f_\lambda))$ converges to $A\in CL(X)$, and in case $A\neq X$, for each $K\in \mathbb{K}(X\setminus A)$, there exists $\lambda_K$ such that $K\subseteq dom(f_{\lambda})$, for each $\lambda > \lambda_K$ and $(f_\lambda|_K)_{\lambda > \lambda_K}$ 
    is uniformly Cauchy. In this case, we will say that the sequence $(f_\lambda|_K)$ is \textit{eventually uniformly Cauchy}.
\end{defn}
\begin{lema}\label{rrr}
    A net $(f_\lambda)_\lambda$ in $C_{od}(X,Y)$ is $\gamma-$Cauchy if, and only if, $(f_\lambda)_\lambda$ is $\tau_{_{\iota,D}}$-converging to some $f\in C_{od}(X,Y)$.
\end{lema}
\begin{proof}
    It is enough to prove the necessary condition. 
    Let $(f_\lambda)_\lambda$ be a $\gamma$-Cauchy net in $C_{od}(X,Y)$. Then $\xymatrix{D(f_\lambda) \ar[r]^-{\tau_{_{Fell}}} & A}$. If $A=X$, $\xymatrix{f_\lambda \ar[r]^-{\tau_{_{\iota,D}}} & \emptyset}$. If $A\neq X$, since $X$ is locally compact and second countable space and $X\setminus A$ is a non-empty open subset of $X$, for each $n\in\mathbb{N}$, there exists $K_n\in \mathbb{K}(X\setminus A)$ such that $$K_n\subseteq int(K_{n+1})\ \text{ and }\ \bigcup\limits_{n\in 
    \mathbb{N}}K_n=X\setminus A.$$
    Since $(f_\lambda)$ is $\gamma$-Cauchy, for each $n\in \mathbb{N}$ there exists $\lambda_n$ such that  $(f_\lambda|_{K_n})_{\lambda>\lambda_n}$ is uniformly Cauchy, then $(f_\lambda|_{K_n})_{\lambda>\lambda_n}$ converges uniformly to $f_n:K_n\longrightarrow Y,$ for each $n\in\mathbb{N}$. Since $f_1\subseteq f_2\subseteq\cdots$ we can define the union function
    $$f=\bigvee\limits_{n\in \mathbb{N}} f_n$$
    Note that $\mathrm{dom}(f)=\bigcup\limits_{n\in \mathbb{N}}K_n=X\setminus A$ (\textit{i.e.} $\xymatrix{D(f_\lambda) \ar[r]^-{\tau_{_{Fell}}} & D(f)}$). If $K\in \mathbb{K}(\mathrm{dom}(f)),\ \epsilon >0$ and fix $n_{_0}\in \mathbb{N}$ such that $K\subseteq K_{n_{_0}}\subseteq dom(f)$. Since $(f_\lambda)$ is $\gamma-$Cauchy, there exists $\lambda_0$ such that $K\subseteq K_{n_0}\subseteq \mathrm{dom}(f_\lambda),$ for each $\lambda >\lambda_0$ and $(f_\lambda|_{K_{n_0}})_{\lambda>\lambda_0}$ is uniformly Cauchy. This net converges to $f_{n_0}=f|_{_{K_{n_0}}}$, \textit{i.e.} there exists $\lambda_1>\lambda_0$ such that $$d_K(f_\lambda,f)\leq d_{K_{n_0}}(f_\lambda,f)<\epsilon,\text{ for all } \lambda>\lambda_1$$
    Hence $f_\lambda\in B_K(f,\epsilon),$ for each $\lambda>\lambda_1$. That is, $\xymatrix{f_\lambda \ar[r]^{\tau_{cc}} & f}$ and therefore $\xymatrix{f_\lambda \ar[r]^-{\tau_{\iota,D}} & f.}$ $f$ is continuous by the Lemma \ref{coro}.
\end{proof}
\begin{thm}
   $(C_{od}(X,Y),\beta)$ is a complete metric space.
\end{thm}
\begin{proof}
    Let $(f_\lambda)$ $\beta$-Cauchy in $C_{od}(X,Y)$, then for the Observation \ref{aaa}, $(D(f_\lambda))$ is $d_{_{Fell}}$-Cauchy and since
    $(CL(X),\tau_{_{Fell}})$ is a compact space, $\xymatrix{D(f_\lambda) \ar[r]^-{\tau_{Fell}}  & A}$. If $A=X$, $\xymatrix{f_\lambda \ar[r]^-{\tau_{\iota,D}}  & \emptyset}$. If $A\neq X$, let $K\in \mathbb{K}(X\setminus A)$ and $\epsilon>0$. There exist $M,N\in \mathbb{N}$ such that $K\subseteq K_{MN}\subseteq U_N\subseteq X\setminus A$. Since $K_{(M+1)N}\subseteq X\setminus A$ and $(f_\lambda)$ is $\beta$-Cauchy there exists $\lambda_0$ such that $K_{(M+1)N}\subseteq D_{f_\lambda}$ and $\beta(f_\lambda,f_\sigma)<\epsilon\, 2^{-(M+N)},$ for each $\lambda,\sigma>\lambda_0$. Hence
    \begin{equation*}
        \begin{split}
            d_K(f_\lambda,f_\sigma) & \leq d_{K_{MN}}(f_\lambda,f_\sigma)\\
            &= \beta_{MN}(f_\lambda,f_\sigma)\\
            &=2^{M+N} 2^{-(M+N)}\beta_{MN}(f_\lambda,f_\sigma)\\
            &\leq2^{M+N}\beta(f_\lambda,f_\sigma)\\
            &<2^{M+N}\epsilon2^{-(M+N)}=\epsilon.
        \end{split}
    \end{equation*}
Therefore $(f_\lambda)$ is $\gamma$-Cauchy and using the lemma \ref{rrr} $\tau_{\iota,D}$-converges.
\end{proof}
As a consequence of the observation \ref{2n} and the previous theorem we have that
\begin{coro}
    If $X$ and $Y$ are locally compact and second contable spaces with $Y$ completelly metrizable space, then $(C_{od}(X,Y),\tau_{\iota,D})$ is a Polish space. 
\end{coro}
Pitifully $(\Gamma [0,1],d_{_\gamma})$ isn't a complete metric space since the 
sequence $f_n(x)=nx,\ x\in\left[0,\frac{1}{n}\right)$ $\tau_{_{\iota,D}}$ (and therefore $\beta$) converges to $\emptyset$, because $\xymatrix{D(f_n) \ar[r]^{\tau_{_{Fell}}} & [0,1].}$ Since that $f_n^{-1}(x)=\frac{1}{n}x,\ x\in [0,1)$ and $\xymatrix{f_n^{-1} \ar[r]^{\beta} & g}$, $g(x)=0,\ x\in[0,1)$. Therefore $(f_n)\subseteq \Gamma[0,1]$ is $d_{_\gamma}$-Cauchy and it doesn't $d_{_\gamma}$-converge.

Although the question of whether $(\Gamma(X),\tau_{hco})$ is a Polish space remains open, the following theorem aims to present an intuition for finding a Polish metric in this very elusive space.
\begin{thm}
    Let $(f_\lambda)$ be a sequence in $\Gamma(X)$ such that  $\xymatrix{f_\lambda \ar[r]^{\tau_{cc}} & f}$, $\xymatrix{f_\lambda^{-1} \ar[r]^{\tau_{cc}} & g}$, $f,g\in C_{od}(X,X)$, with  $\mathrm{im}(f)\subseteq \mathrm{dom}(g)$ and $\mathrm{im}(g)\subseteq \mathrm{dom}(f)$. Then $f$ and $g$ are homeomorphisms and $f^{-1}=g$.  
\end{thm}
\begin{proof}
    Let $(f_\lambda)$ such sequence, $x\in \mathrm{dom}(f)$ and $\epsilon>0$. There exists $\lambda_1$ such that $x\in dom(f_\lambda)$ and $(f_\lambda(x))_{\lambda >\lambda_1}$ converges to $f(x)\in \mathrm{im}(f)\subseteq \mathrm{dom}(g)$. Since $\mathrm{dom}(g)$ is open subset of $Y$, there exists $\lambda_2>\lambda_1$ such that $f_{\lambda}(x)\in \mathrm{dom}(g)$, for each $\lambda>\lambda_2$. Given that $L=\{f_\lambda(x) ,f(x):\lambda>\lambda_2\}$ is a compact subset of $\mathrm{dom}(g)$ and $\xymatrix{f_\lambda^{-1} \ar[r]^{\tau_{cc}} & g}$ there exists $\lambda_3>\lambda_2$ such that $L\subseteq dom(f_\lambda^{-1})=im(f_\lambda)$ and 
    \begin{equation*}
        d_L(f_\lambda^{-1},g)<\frac{\epsilon}{2}, \text{ for each } \lambda>\lambda_3.
    \end{equation*}
    Given that $g$ is continuous and $(f_\lambda(x))_{\lambda>\lambda_3}$ converges to $f(x)\in \mathrm{dom}(g)$ there exists $\lambda_0>\lambda_3$ such that 
    \begin{equation*}
        d(g(f_\lambda(x)),g(f(x)))<\frac{\epsilon}{2}, \text{ for each }\lambda>\lambda_0.
    \end{equation*}
    So, for each $\lambda>\lambda_0$ we have that $$d(x,g(f(x)))\leq d(f_\lambda^{-1}(f_\lambda(x)),g(f_\lambda(x)))+d(g(f_\lambda(x)),g(f(x)))<\frac{\epsilon}{2}+\frac{\epsilon}{2}=\epsilon$$
    therefore, $g(f(x))=x,$ for each $x\in \mathrm{dom}(f)$. In an analogous way we can show that $f(g(y))=y,$ for each $y\in \mathrm{dom}(g)$. How $f$ and $g$ are continuous functions we conclude the result.
\end{proof}

Note that the function $\varphi:\Gamma(X)\longrightarrow C_{od}(X)^2,\ f\longmapsto (f,f^{-1})$ is an $\tau_{_{hco}}$-$\,\tau_{_{\iota,D}}\times\tau_{_{\iota,D}}$ embedding. In particular $\varphi(\Gamma[0,1])$ is neither open nor closed subset in $(C_{od}[0,1]^2,\tau_{_{\iota,D}}\times\tau_{_{\iota,D}})$.
\begin{center}
    \textit{\textbf{Conjeture.} $\varphi$ is an $G_\delta$ embedding.} 
\end{center}

\end{document}